\documentclass[12pt,reqno]{amsart}
\usepackage{amssymb,color}
\usepackage{graphicx}

\newcommand{\e}{\varepsilon}

\newcommand{\la}{\lambda}
\newcommand{\al}{\alpha}

\newcommand{\fy}{\varphi}

\newcommand{\p}{\partial}
\newcommand{\I}{\infty}

\newcommand{\R}{\mathbb{R}}

\renewcommand{\S}{\mathcal{S}}

\numberwithin{equation}{section}

\newtheorem{thm}{Theorem}[section]

\newtheorem{lem}[thm]{Lemma}

\theoremstyle{remark}

\newtheorem{defn}{Definition}

\newcommand{\ran}{\rangle}
\newcommand{\lan}{\langle}

\newcommand{\EQ}[1]{\begin{equation} \begin{split} #1 \end{split} \end{equation}}
\newcommand{\Del}[1]{}

\newcommand{\D}{\mathcal{D}}

\newcommand{\HH}{\mathcal{H}}

\def\PS{\mathcal{PS}}

\def\dim{\mathrm{dim}}
\def\ST{\mathrm{ST}}
\def\PS{\mathrm{PS}}
\def\nn{\nonumber}
\def\f{\frac}
\def\lan{\langle}
\def\ran{\rangle}
\def\exp{\mathrm{exp}}
\def\const{\mathrm{const}}

\begin{document}

\author{R.~Donninger}
\address{Department of Mathematics, The University of Chicago\\ Chicago, IL 60615, U.S.A.}
\email{donninger@uchicago.edu}

\author{W.~Schlag}
\address{Department of  Mathematics, The University of Chicago\\ Chicago, IL 60615, U.S.A.} 
\email{schlag@math.uchicago.edu} 

\thanks{The first author was supported by a Schr\"odinger Stipendium by the Austrian FWF, Project No.~J2843. The second author was supported in part by the National Science Foundation,  DMS-0617854  as well as by a Guggenheim fellowship.  This research was supported in part by the NSF through {\em Teragrid}
resources provided by NICS, the {\em National Institute for Computational Sciences} under grant number TG-DMS110003.
The authors thank  Piotr Bizo\'n, Dugald Duncan, Frank Merle,  and Ridgway Scott for discussions and useful advice. }

\title[Numerical Study of  NLKG]{Numerical Study of the blowup/global existence dichotomy for the focusing cubic nonlinear Klein-Gordon equation}


\begin{abstract}
We present some numerical findings concerning the nature of the  blowup vs.\ global existence dichotomy for the focusing cubic nonlinear
Klein-Gordon equation in three dimensions for radial data. The context of this study is provided by the classical paper by 
 Payne, Sattinger~\cite{PS}, as well as the recent work by K.~Nakanishi, and the second author~\cite{NakS1}. 
Specifically, we numerically investigate the boundary of the forward scattering region. 
While the results of \cite{NakS1} guarantee that this boundary is smooth at energies which are near the ground state energy, 
it is currently unknown whether or not it continues to be a smooth manifold at higher energies. 
While we do not find convincing evidence of either smoothness or singularity formation, our numerical work does indicate that 
at larger energies the boundary becomes much more complicated than 
at energies near that of the ground state. 
\end{abstract}

\maketitle


\section{Introduction}

This paper is concerned with  the nonlinear Klein-Gordon equation 
\EQ{
\label{eq:NLKG3}
\Box u + u = u_{tt}-\Delta u+u = u^{3}, \quad (t,x)\in \R^{1+3}
}
with radial data.    
This equation exhibits several conserved quantities, of which the  energy, with $\vec u=(u,\dot u)$, 
\EQ{\label{eq:energy}
E(\vec u) = \int_{\R^{3}} \big[\frac12( |\nabla u|^{2}+u^{2}+\dot u^{2})-\frac{1}{4}|u|^{4}\big]\, dx 
}
is the most important one for our purposes.  Note that we are only considering the focusing case here. 
Let us briefly review some basic well-known facts concerning the well-posedness of~\eqref{eq:NLKG3} in
the energy space (this does not require any radial assumption). 

\begin{thm} \label{thm:basic}
The NLKG equation~\eqref{eq:NLKG3}   is locally wellposed for data in $H^{1}(\R^{3})\times L^{2}(\R^{3})$ in the usual (Duhamel
formulation) sense and the energy is conserved. 
Small data lead to global existence and scattering to zero, whereas data of negative energy lead to
finite time blowup in either time direction.  If an energy solution exists on $[0,T_{*})$ with $T_{*}<\I$ maximal, then necessarily
$\| u\|_{L_t^{3}([0,T_{*}); L_x^{6} )  }=\I$; if $T_{*}=\I$ and $\|u\|_{L^{3}_{t}([0,\I);L^{6}_{x})}<\I$, then $u$ scatters
to zero as $t\to \I$. Finally, smooth data  lead to smooth solutions. 
\end{thm}

The proof  is a straightforward application of energy and Strichartz estimates, see for example~\cite{Strauss}.  In the defocusing case
one has global existence and scattering to zero for all data, see~\cite{Jorgens}, \cite{Bren1}, \cite{Bren2}, \cite{GV1}, \cite{GV2}, \cite{Pech}. 

In what follows, $\HH:=H^{1}(\R^{3})\times L^{2}(\R^{3})$ is the phase (energy) space, and $\|\cdot\|_{\ST}$ stands for
the $L^{3}_{t} L^{6}_{x}$-Strichartz norm associated with~\eqref{eq:NLKG3}, respectively, on the time-interval $[0,\I)$. 
Let us now define  the forward scattering set.

\begin{defn}
The forward scattering set is defined as 
\EQ{\label{eq:S+def}
\S_{+}:= \Big\{ (u_{0}, u_{1})\in\HH\mid u(t):=S(t) (u_{0}, u_{1}) \; \exists \;\forall t\ge0,   \| u\|_{\ST} < \I \Big\} 
}
where $S(t)$ is the nonlinear evolution of the NLKG equation.
\end{defn}

The following properties of $\S_{+}$ follow from 
 the main wellposedness theorem above and some simple perturbative arguments.

\begin{lem}
The set $\S_{+}$ enjoys the following properties:
\begin{itemize}
\item $\S_{+}\supset B_{\delta}(0)$, a small ball in~$\HH$
\item  $\S_{+}$ is an open set
\item  $\S_{+}\neq\HH$
\item $\S_{+}$ is path-connected. 
\end{itemize}
\end{lem}

For the third property one can use that solutions of negative energy blowup in finite time, see~\cite{Levine}, \cite{PS}. 
The following questions pose themselves quite naturally now:

\begin{enumerate}

\item Is $\S_+$ bounded in $\HH$?

\item  What does $\p\S_{+}$  look like, is  it a smooth manifold, or can it be very rough?

\item If $\p\S_{+}$ is a smooth manifold, does it separate a region of finite time blowup from one of global existence,
at least locally? 

\item What is the dynamical behavior of  solutions starting on $\p\S_{+}$? Are there any special solutions with data on $\p\S_{+}$?

\end{enumerate}

Some partial answers to these questions were found in~\cite{NakS1} for~\eqref{eq:NLKG3} 
with radial data (the $\dim=3$ case with nonradial data was treated in~\cite{NakS3}, the cubic NLS equation in $\dim=3$ in~\cite{NakS2}, and
the critical wave equation in $\dim=3,5$ in~\cite{KNS1}, whereas the one-dimensional case is studied in~\cite{KNS2}). 
To be more specific, it was found that the answer to (1) is ``no''. In fact,
there exists a curve in $\HH$ connecting $0$ to~$\I$ and such that a solution starting from an arbitrary point on that curve belongs to $\S_{+}$. 
The answers to the other questions are formulated in terms of the {\em ground state} $Q$. This refers to the unique positive, radial, stationary 
energy solution of the NLKG equation. In other words, $Q=Q(r)>0$, $Q\in H^{1}(\R^{d})$, and 
\EQ{\label{Qeq}
-\Delta Q + Q = Q^{p}
}
Amongst all solutions $\fy$ of~\eqref{Qeq}, the ground state minimizes the energy
\[
J(\fy):= \int_{\R^{d}} \big(\frac12[|\nabla \fy|^{2}+\fy^{2}]-\frac{1}{p+1} |\fy |^{p+1}\big)\, dx
\] 
In $\dim=1$, \eqref{Qeq} has only two solutions (up to translation symmetry) which decay at $\pm\I$, namely
\EQ{
\label{eq:soliton}
Q(x) =  \alpha \cosh^{-\frac{1}{\beta}}(\beta x),\qquad \al=\big(\frac{p+1}{2}\big)^{\frac{1}{p-1}}, \;\; \beta=\frac{p-1}{2}
}
In contrast, in $\dim=3$, \eqref{Qeq} has infinitely many $H^{1}$ solutions, but up to translation symmetry only one positive solution, namely
the ground state, see~\cite{Strauss77}, \cite{BerLions}, \cite{Coff}. 
In their classical paper~\cite{PS}, Payne and Sattinger showed the following, with
\[
K(\fy):=\p_{\la} J(e^{\la}\fy)\Big|_{\la=0} = \int_{\R^{d}}\big[ |\nabla \fy|^{2} + |\fy|^{2} - |\fy|^{p+1} \big](x)\, dx
\]
being the scaling functional. 

\begin{thm}
\label{thm:trap} The two regions
\EQ{\nn
\PS_{+} := \{ (u_{0},u_{1})\in\HH  \mid E(u_{0}, u_{1}) < J(Q),\; K(u_{0})\ge 0\} \\
\PS_{-}  := \{ (u_{0},u_{1})\in\HH  \mid E(u_{0}, u_{1})   < J(Q),\; K(u_{0}) < 0\}
}
are invariant under the flow of~\eqref{eq:NLKG3} in the following sense: if $(u(0),\dot u(0))\in \PS_{+}$, then $(u(t),\dot u(t))\in \PS_{+}$
for as long as the solution exists, and the same holds for $\PS_{-}$.  Solutions of \eqref{eq:NLKG3} which lie in~$\PS_{+}$  
exist for all times, whereas those in $\PS_{-}$ blow up in finite time (in both temporal directions).  In particular, data of negative energy blow up in finite 
time, and $\S_{+}\neq\HH$. 
\end{thm}

The invariance is a fairly immediate consequence of the fact that  the variational problem
\[
\inf\{ J(\fy) \mid \fy\in H^{1}\setminus\{0\}, K(\fy)=0\}
\]
has $\pm Q$ as the unique solutions (up to translation). Furthermore, in $\PS_{+}$ one has 
\[
E(\vec u) \simeq \| \vec u\|_{\HH}
\]
whence Theorem~\ref{thm:basic} implies global existence. The blowup in $\PS_{-}$ follows from a concavity argument. 
Very recently,  Ibrahim, Masmoudi, Nakanishi~\cite{IMN} proved that solutions in $\PS_{+}$ scatter to zero as $t\to\pm\I$. Thus, $\PS_{+}\subset \S_{+}$. 
Returning to describing the partial answers to questions (2), (3), (4) above, introduce
\[
\HH^{\e}:= \{(u_{0},u_{1})\in \HH\mid E(u_{0},u_{1})<J(Q)+\e^{2}\}
\]
Then the following was shown in \cite{NakS1}, \cite{KNS2}:

\begin{thm}\label{thm:NS}
 For some small $\e>0$ one has 
 \begin{itemize}
 \item  $\p\S_{+}\cap \HH^{\e}$ is a smooth, co-dimension one
manifold
\item it does separate a region of scattering to zero as $t\to\I$ from one exhibiting finite time blowup
\item any data 
from $\p\S_{+}\cap \HH^{\e}$  lead to global solutions in forward time which scatter to $\pm Q$, i.e., 
\EQ{\label{eq:scattoQ}
\vec u = \pm (Q,0) + \vec v + o_{\HH}(1) \qquad t\to\I
}
where $\vec v$ is an energy solution of the free Klein-Gordon equation. 
\end{itemize}
\end{thm}

These results are for radial data, see~\cite{NakS3} for a version
with nonradial data. Furthermore, the aforementioned references contain a description of the boundary $\p\S_{+}\cap \HH^{\e}$ as
the center-stable manifold associated with $(Q,0)$ in the sense of Bates, Jones~\cite{BJ}. Moreover, the stable (and unstable) manifolds are
one-dimensional and are described in the terms of the threshold solutions found by Duyckaerts, Merle in the energy critical case~\cite{DM1}, \cite{DM2}. 

A well-known consequence of Theorem~\ref{thm:trap} is the instability of $Q$. More precisely, in any neighborhood of $(Q,0)$ in $\HH$ we
can find data which lead to finite time blowup as well as global existence, respectively. To see this, we remark that the linearized 
operator $L_{+}=-\Delta+1-pQ^{p-1}$ has a single negative eigenvalue, the ground state. It is evident that $L_{+}$ has negative spectrum since
\[
\lan L_{+} Q|Q\ran <0.
\] 
Let $L_{+}\rho=-k^{2}\rho$ where $k\ne0$ and $\| \rho \|_{2}=1$. 
For further discussion of the spectral properties of $L_{+}$ and $L_{-}=-\Delta+1-Q^{p-1}$, see \cite{KrS1} for the one-dimensional case, as well as~\cite[Lemma 2.3]{NakS1}, \cite{DS}, \cite{CHS}
 for the three-dimensional case.  Expanding the stationary energy $J$ and the functional $K$ around $Q$ yields
\EQ{
J(Q+v) &= J(Q) + \frac12 \lan L_{+} v| v\ran + O(\| v\|_{H^{1}}^{3})  \\
K(Q+v) &=  -(p-1) \lan Q^{p}| v\ran + O(\| v\|_{H^{1}}^{2})
}
as $\|v\|_{H^{1}}\to0$. Therefore, by Theorem~\ref{thm:trap}, data $(Q+\e \rho,0) $ blow up for $\e>0$ small and scatter to zero for $\e<0$ small. 
The same holds for $(1+\e)Q$. 

Expanding further, we set $u=Q+\la \rho + w$ where $w\perp \rho$. Then  $\dot u = \dot\la \rho + \dot w$ and 
\EQ{
E(\vec u) = J(Q)  + \frac12(\dot\la^{2}-k^{2}\la^{2}) + \frac12 (\lan L_{+}w|w\ran + \|\dot w\|_{2}^{2}) + O(\|v\|_{H^{1}}^{3}) 
}
where $v=\la \rho+w$. Since $\lan L_{+}w|w\ran\ge0$, we conclude that the only way to decrease the energy below $J(Q)$ is through $-k^{2}\la^{2}$.
Moreover, in a $\delta$-neighborhood of $(0,0)$ in the $(\la,\dot\la)$-plane the set $\{E(\vec u)<J(Q)\}$ looks like 
\EQ{
\label{eq:cones}
\{ (\xi,\eta)\in\R^{2}\mid \xi^{2}-\eta^{2}<0,\; |\xi|+|\eta|<\delta\}
}
at least up to cubic corrections. 

Needless to say that there are other important aspects of the problem which we do not touch upon 
in the present work. 
For instance, it is known \cite{HZ} that the blowup rate is
determined by the ODE $u_{tt}=u^3$, very similar to the analogous situation for the wave equation \cite{MZ}.

\section{The numerical methods}
\label{sec:code}

At the core of all numerical work lie  second order finite difference schemes (both explicit and implicit schemes were used) which
solve the cubic NLKG equation with radial data.  Such schemes are quite standard and are discussed in more detail 
below. 
In order to study the set $\S_{+}$ for~\eqref{eq:NLKG3} numerically, the authors applied  these solvers to 
data of the form $(Q+Af, Bg)$ or $(Af, Bg)$ where $f,g$ are fixed radial functions, such as Gaussians, exponentials,
with or without oscillatory factors (typical ones are $\sin(6r)$ or $\sin(6r^{2})$). Here $(A,B)$ runs over a fine grid of points from a suitable rectangle. 
If numerical blowup is found, which simply means that the $H^{1}\times L^{2}$ norm on some fixed  ball becomes very large in finite time, something
like $10^{5}$ times larger than the norm of the initial data, the pair $(A,B)$ is stored in a blowup file;  if the $H^{1}\times L^{2}$-norm becomes
very small on some ball (say, falls below $1/10$ of the size of the data), and remains small, then the pair $(A,B)$ is regarded as ``dispersive''. If no decision of this type can
be reached up to some prescribed maximum number of time steps, the data set is characterized as ``indecisive''.  Note that the latter
case occurs if $(Q+Af, Bg)$ falls exactly on the center-stable manifold (an extremely unlikely event; however, one may come very near to this manifold,
which then requires a large number of time-steps in order to reach the blowup/dispersion decision).   

While this characterization is somewhat delicate (for example, it is not apriori clear why solutions could not blow up on a sphere of 
large radius, cf.~\cite{Raphael}), the authors have found empirically that it leads to consistent and meaningful results (more on this can be
found in Section~\ref{sec:resultsd3}). 

The numerical calculations were performed with difference schemes applied to the equation
\EQ{
\label{eq:radNLKG3}
v_{tt} - v_{rr} + v = \frac{1}{r^{2}} v^{3}
}
for $v(t,r)=r u(t,r)$.  The choice of the difference scheme is a somewhat delicate matter. For the free wave equation there are several   schemes,
accurate to second or higher orders, which are known to converge~\cite{Strik}. For the nonlinear case, 
especially with the singularity on the right-hand side of~\eqref{eq:radNLKG3}, there does not appear to be much known in terms of the convergence of difference schemes for~\eqref{eq:radNLKG3}. 

The first, implicit, second order accurate difference scheme which we use is due to Strauss, Vazquez~\cite{StrV} (which they introduced for the defocusing equation) and reads as follows: 
\EQ{\label{eq:StrV}
&\frac{v_{j}^{n+1} -2 v_{j}^{n} + v_{j}^{n-1}}{(\Delta t)^{2}} - \frac{v_{j+1}^{n} - 2v_{j}^{n} + v_{j-1}^{n}}{(\Delta r)^{2}} +\frac12[v_{j}^{n+1}+ v_{j}^{n-1}] \\
&\quad  - \frac{1}{4(j\Delta r)^{2}} \big[ (v_{j}^{n+1})^{3} + (v_{j}^{n+1})^{2} v_{j}^{n-1} + v_{j}^{n+1} (v_{j}^{n-1})^{2}
+(v_{j}^{n-1})^{3}\big] =0
}
the nonlinear term being $\frac{G(v_{j}^{n+1})- G(v_{j}^{n-1})}{v_{j}^{n+1} - v_{j}^{n-1}}$ with $G(v)=-\frac14 v^{4}$. 
Here, we use the standard shorthand notation $v^n_j:=v(n\Delta t, j\Delta r)$.
The discretization steps $\Delta t, \Delta r$
are chosen such that $\frac{\Delta t}{\Delta r}=0.9$; in particular, the Courant-Friedichs-Lewy condition is satisfied. Moreover, typical values of $\Delta r$ are on the order of $10^{-3}$. 
For the implementation, one applies a Newton scheme to solve for $v_{j}^{n+1}$ using as a starting guess the explicit version of~\eqref{eq:StrV} in
which all but the first occurrences of $v_{j}^{n+1}$ are replaced by $v_{j}^{n}$, say. Note that due to the  focusing sign the Newton scheme
may become degenerate, which means that the derivative of the polynomial of $v_{j}^{n+1}$ as defined by~\eqref{eq:StrV} can vanish or almost vanish. 
However, one checks that this can only happen if $M(t)\cdot (\Delta r)>1$ where $$M(t)=\max_{r>0} |u(t,r)|=\max_{r>0}|v(t,r)|/r;$$ if the latter case occurs, then 
we terminate with a decision of blowup. 
To calculate the solution of~\eqref{eq:radNLKG3} on a rectangle $(t,r)\in (0,T)\times (0,R)$ we use a space-time grid on the larger rectangle $(0,T)\times (0,R+T)$,
since the latter is the smallest rectangle containing the domain of dependence of~$(0,T)\times (0,R)$. 

The scheme \eqref{eq:StrV}  conserves the discrete energy
\EQ{\label{eq:Ediscr}
E_{n} &:= \Delta r\cdot\Big[ \frac12\sum_{j} \Big(\frac{v_{j}^{n+1}-v_{j}^{n}}{\Delta t}\Big)^{2} + \frac12 \sum_{j} \Big(\frac{v_{j+1}^{n+1}-v_{j}^{n+1}}{\Delta r}\Big)
\Big(\frac{v_{j+1}^{n}-v_{j}^{n}}{\Delta r}\Big) \\
&\qquad +\frac12 \sum_{j} \frac{( v_{j}^{n+1})^{2} +  ( v_{j}^{n})^{2} }{2} - \sum_{j} \frac{ (v_{j}^{n+1})^{4} + (v_{j}^{n})^{4} }{8 (j \Delta r)^{2} } 
\Big]
}
and is stable and thus convergent~\cite{Strik} for the free equation. 

The second, explicit, second order accurate difference scheme is as follows:
\EQ{\label{eq:Rol}
\frac{v_{j}^{n+1} -2 v_{j}^{n} + v_{j}^{n-1}}{(\Delta t)^{2}} - \frac{v_{j+1}^{n} - 2v_{j}^{n} + v_{j-1}^{n}}{(\Delta r)^{2}} +v_{j}^{n} 
 - \frac{1}{(j\Delta r)^{2}} (v_{j}^{n})^{3} =0
}
It is easier to implement, and runs faster since no Newton iteration is required. It is
still stable and convergent for the free equation.  While it does not conserve the discrete energy~\eqref{eq:Ediscr} {\em exactly},  the authors have found that 
it does so {\em approximately} in those examples which they considered.  

A theoretical comparison of the difference schemes~\eqref{eq:StrV} and~\eqref{eq:Rol} appears to be quite challenging. 
However, after running numerous 
computations with these schemes the authors have found that they lead to the same conclusions as far
as the dispersion/blowup dichotomy is concerned (for example, in terms of the  pictures appearing in the following section). 
However, as far as the actual blowup solutions are concerned, significant differences can be found near the blowup time. But we reiterate that in   all cases checked by the authors 
the qualitative behavior comes out the same, which means that the  appearance of the sets in Section~\ref{sec:resultsd3} comes out the same. 

Other than running their computations with both \eqref{eq:StrV} and~\eqref{eq:Rol} (the latter typically with finer resolution due to its faster execution), and then repeating
the calculations with $\Delta r/2$, 
the authors have used the {\em bisection method} as a means of checking the validity of figures such as those appearing in Section~\ref{sec:resultsd3}.
This refers to a process in which two data pairs $d_{1}:=(A_{1},B_{1})$ and $d_{2}:=(A_{2}, B_{2})$, the first leading to blowup, and the second to dispersion, are chosen in
close proximity to each other. Then one bisects the line segment joining them, tests the midpoint for blowup/dispersion, and then   continues with the line segment connecting
the midpoint to either $d_{1}$ or $d_{2}$, so that the endpoints of this new segment exhibit opposite behaviors. 
Iterating this procedure until machine precision is reached, the authors have found that one obtains data pairs which 
``almost'' lie on the center-stable manifold. This means that the time evolution of these data exhibit ``metastable'' behavior, or ``quasinormal'' ringing.
The latter term is used in  the theoretical physics literature which concerns itself with the decay of dispersive tails, as they appear for example in the context
of the so-called Price law. See also  
Bizo\'n et al.~\cite{Biz}, \cite{BCS}. 

Finally, as an additional check, the authors have also run computations with a code that uses a compactified radial coordinate. 
This approach is quite common in the physics literature and it allows one (at least in principle) 
to follow the evolution up to spatial infinity which is important to rule out 
the possibility of late time blowup. 
Fortunately, such an unpleasant phenomenon does not seem to occur and the Klein-Gordon equation under consideration is well-behaved in this respect.

As far as the coding is concerned, the authors used the ``C''-language. The computations leading to the results of the following section
were carried out on two machines: (1)  a MacPro capable of running 16 processes in parallel. In that case, the authors 
relied on the ``{\tt fork()}''-command to distribute the $(A,B)$-loop over 50 (typically) parallel processes
at any given moment.
(2) The {\em Kraken} supercomputer (a Cray XT5) at NICS, which is part of the {\em TeraGrid}.  

\section{The numerical findings}
\label{sec:resultsd3}

In order to explore questions (2)-(4) above numerically, the authors investigated solutions generated by data sets 
\EQ{ \nn 
\D_{1}(f,g)&:=\{ (Q+Af, Bg)\mid |A|<a, \; |B|<b \} \\
\D_{2}(f,g)&:=\{ (Af, Bg)\mid |A|<a, \; |B|<b \} 
}
where $(f,g)$ are some simple fixed radial rapidly decaying functions, such as Gaussians, exponentials, possibly multiplied with oscillatory factors,
and $a,b>0$ are chosen on a case to case basis. 
\begin{figure}[ht]
\begin{center}
\includegraphics[width=0.496\textwidth]{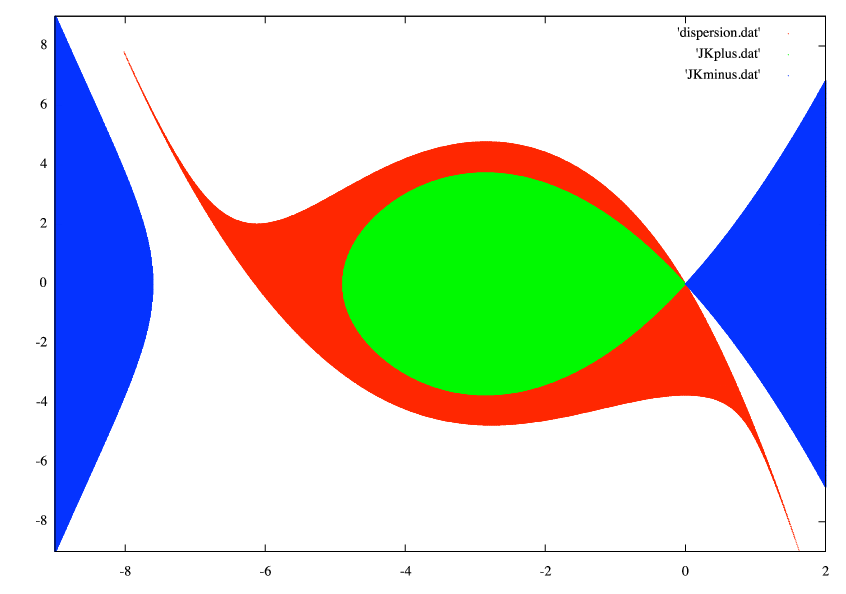}
\includegraphics[width=0.496\textwidth]{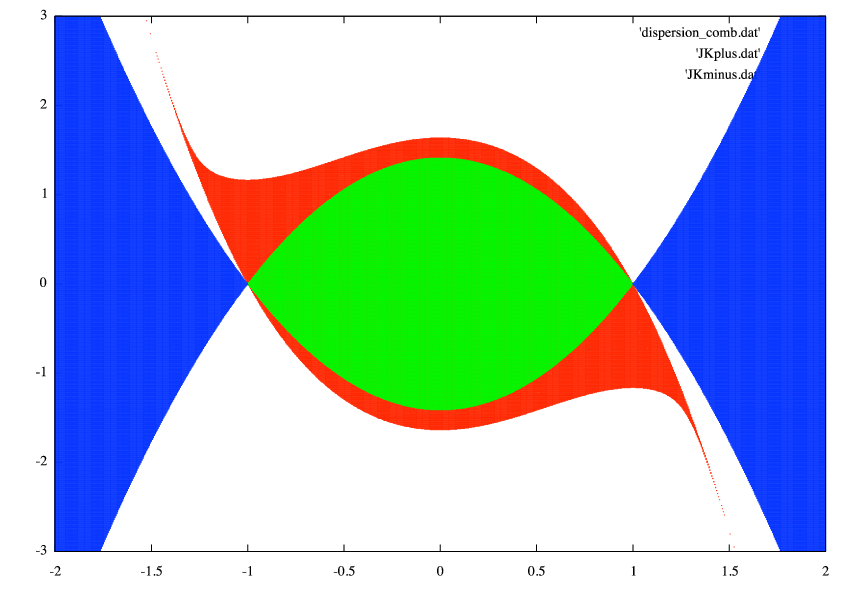}
\end{center}
 \caption{\label{fig:1.1} Numerical results for $\D_{1}(e^{-r^2}, e^{-r^2})$ and $\D_2(Q,Q)$}
\end{figure}
In the examples shown below, $f,g$ are smooth, including near the origin (in fact, they are analytic). 
Once a choice of $(f,g)$ is made, the numbers, $a$, $b$ are then determined by trial and error (weighing numerical cost against 
the desire to obtain a large enough section of $\S_{+}$ and $\p\S_{+}$).  

Discretizing $(A,B)$ by a sufficiently fine rectangular grid, and solving
\eqref{eq:NLKG3} at each grid-point by the difference schemes of the previous section, one then colors each grid-point according to whether 
the numerical solution appears to exhibit blowup (it is kept blank) or scattering to zero (colored red) for large times\footnote{The 
figures shown below represent only a small portion of all data which the authors have produced. The selection shown
here is intended to capture the most important features which can be seen from those data.}.

\begin{figure}[ht]
\begin{center}
 \includegraphics[width=0.496\textwidth]{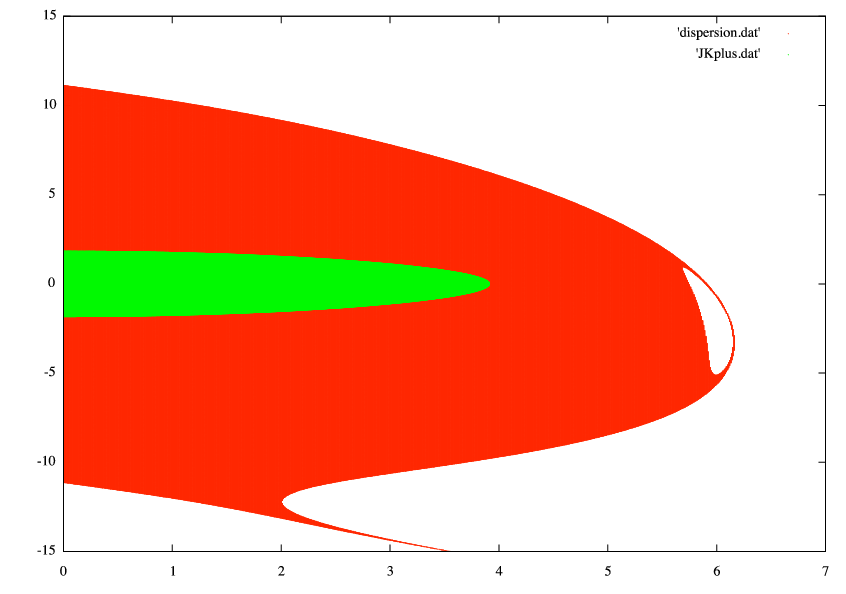}
 \includegraphics[width=0.496\textwidth]{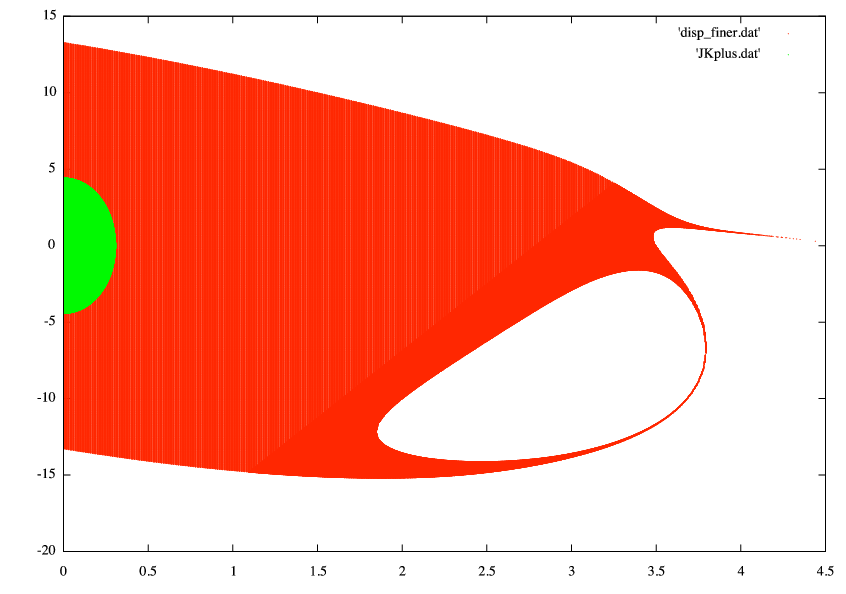}
\end{center}
\caption{\label{fig:1.3}  Numerical results for  $\D_{2} (e^{-\lan r\ran}, \sin(6r^2)e^{-(r^2-1)^2/10})$ and $\D_2(r^2\cos(10r^2)e^{-\lan r^2\ran}, \sin(6r^2)e^{-4(r^2-1)^2})$}
\end{figure}

In this fashion, one obtains two-dimensional sections of $\S_{+}$. In addition, the authors super-imposed the Payne, Sattinger regions $\PS_{\pm}$
over the resulting images so as to obtain a meaningful comparison to Theorem~\ref{thm:trap}.  The image on the left of 
Figure~\ref{fig:1.1} shows the outcome of such a computation for
the data choice 
\EQ{\nn
(u(0),\dot u(0))(r) = (Q(r) + A e^{-r^2}, B e^{-r^2})
}
with the horizontal axis being $A$, and the vertical being $B$. The green and blue regions are $\PS_+$ and $\PS_-$, respectively. Note that
they meet at $(0,0)$, which is precisely the soliton~$Q$, in the conic fashion described by~\eqref{eq:cones}. 
The red region are data that lead to global existence (it includes green as a subset), whereas those that appear white 
(which include blue as a subset) correspond to finite time blowup.  
The image on the right-hand side of Figure~\ref{fig:1.1}  is produced by data $(AQ(r),BQ(r))$. 
Figure~\ref{fig:1.1} agrees with the results of~\cite{PS} and~\cite{NakS1}. 
Indeed,
$\PS_{+}$ does appear as a subset of the red region, $\PS_{-}$ as a subset of the white set, and the boundary of the red region near $(Q,0)$
looks like a smooth curve. It continues to look like a smooth curve even far away from $(Q,0)$, but this is currently not backed by theory.  
While the computation of the green and blue Payne-Sattinger regions does not pose a serious obstacle (the energy $J(Q)$ can by 
obtained once $Q$ is determined by a shooting-method, say, using {\tt Mathematica}), determining the red regions is of course the main
difficulty. For that we use the difference schemes described in Section~\ref{sec:code}. 
\begin{figure}[ht]
\begin{center}
\includegraphics[width=0.496\textwidth]{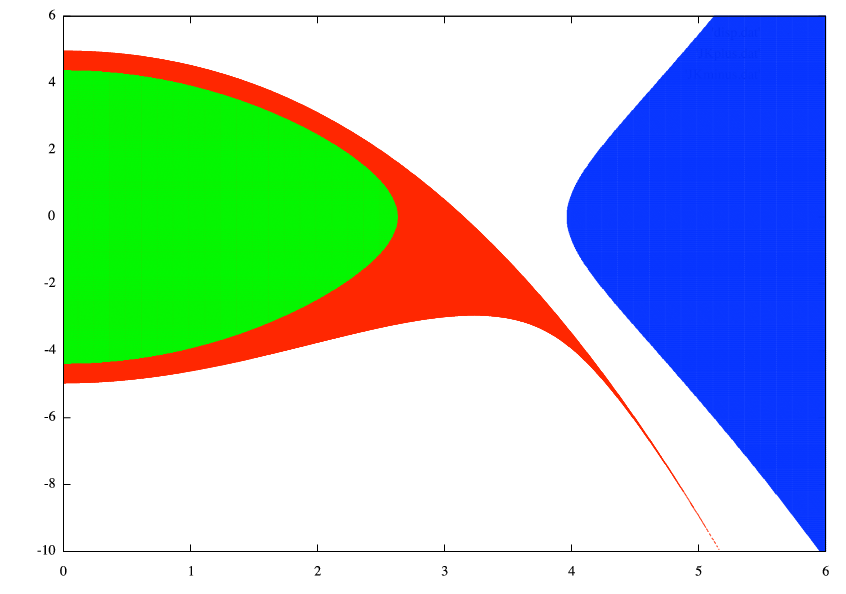}
\includegraphics[width=0.496\textwidth]{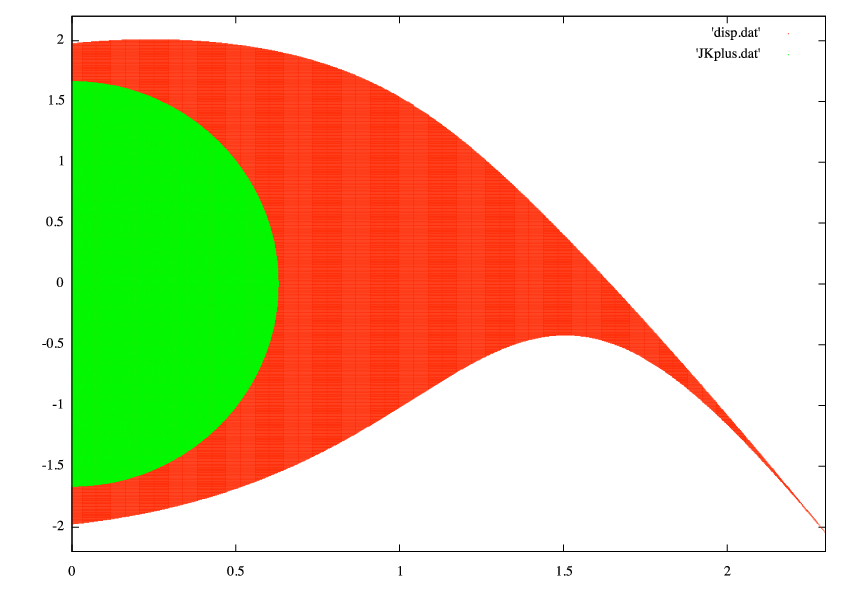}
\end{center}
 \caption{\label{fig:1.2} Numerical results for $\D_{2}(e^{-r^2}, e^{-r^2})$ and $\D_2( r^2 e^{-4(r^2-1)^2} ,  r^2 e^{-(r^2-1)} )$}
\end{figure}
Figure~\ref{fig:1.2} and \ref{fig:1.3} show the outcome of computations where the data do not pass through either $(Q,0)$ or $(-Q,0)$. 
A natural section of this type is  $\D_{2}(e^{-r^{2}}, e^{-r^{2}})$, and it appears on the left-hand side of Figure~\ref{fig:1.2}. The right-hand
image of that figure has very similar features. In particular, "spikes" and a smooth boundary (we discuss these spikes in more details
below). 
Note that in these figures we only consider positive values of~$A$, which suffices by the reflection symmetry about the origin. The blue $\PS_{-}$
region does not intersect the rectangle chosen in the right-hand image of Figure~\ref{fig:1.2}, and thus does not appear in there. 

A common feature of Figures~\ref{fig:1.1}, \ref{fig:1.2}  is that the boundary of the red region appears as a union of smooth curves. 
A new feature in the form of ``bubbles'' is shown in Figure~\ref{fig:1.3}. At least judging from the resolution used in  these images, 
they do not seem to  involve any
singularities. We use the notation $\lan r\ran=\sqrt{1+r^{2}}$ in the definition of the data (we prefer $\lan r\ran$ to $r$ as the former is smooth in $\R^3$,
whereas the latter is not). 

\begin{figure}[ht]
\begin{center}
 \includegraphics[width=0.496\textwidth]{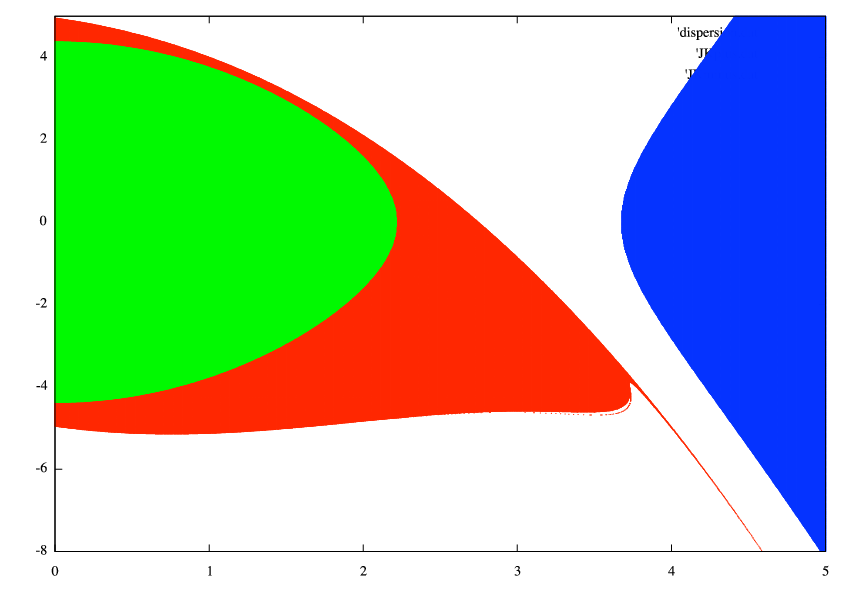}
 \includegraphics[width=0.496\textwidth]{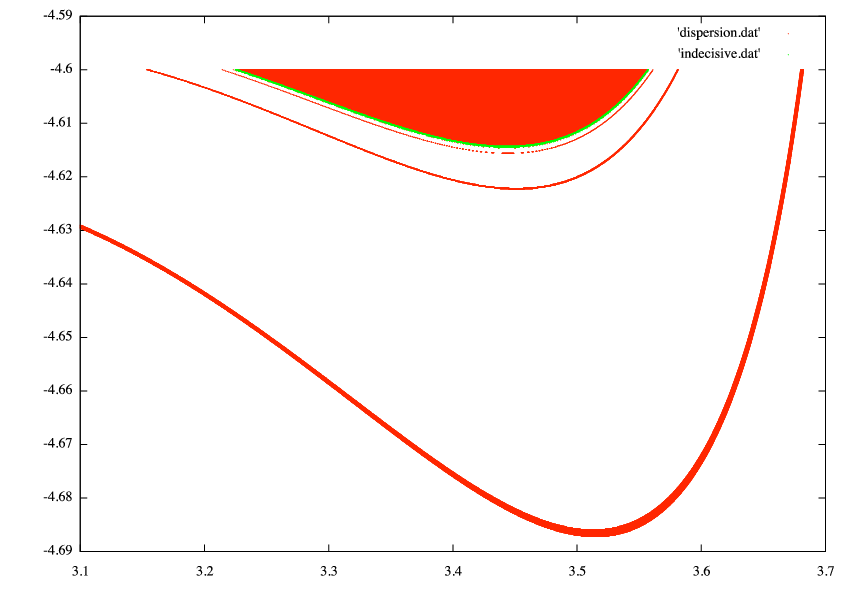}
\end{center}
\caption{\label{fig:1.4} Numerical results for  $\D_{2} (  e^{-\frac{r^2}{\lan r\ran}}, e^{-r^2}  )$ and one zoom near the boundary}
\end{figure}

\begin{figure}[ht]
\begin{center}
 \includegraphics[width=0.495\textwidth]{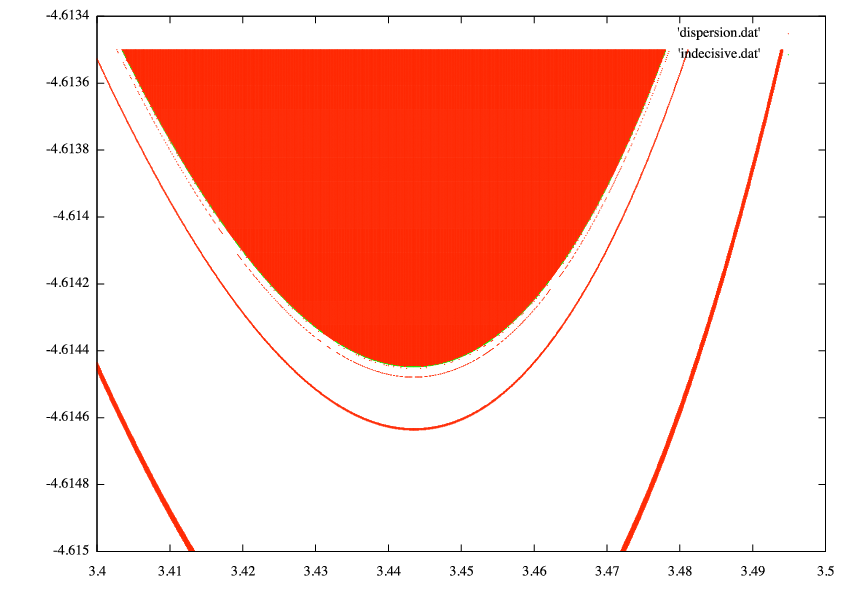}
 \includegraphics[width=0.495\textwidth]{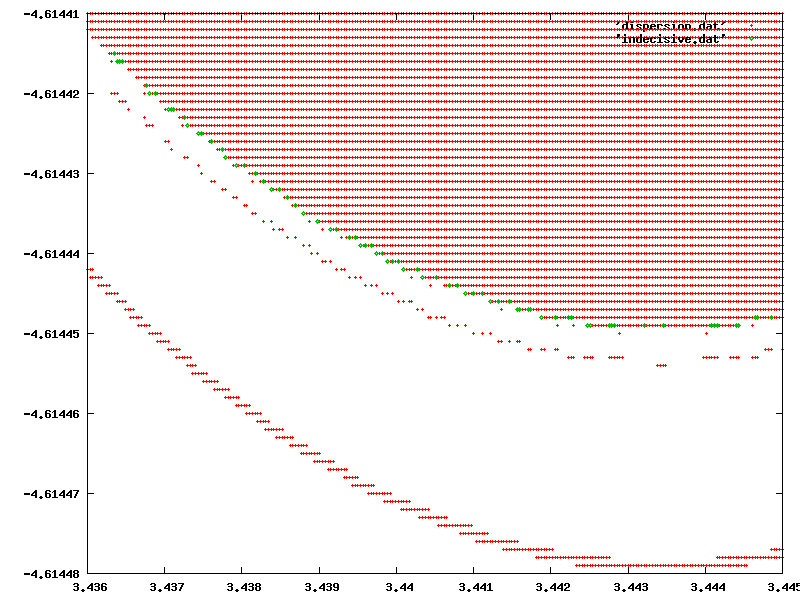}
\end{center}
\caption{\label{fig:1.5}  Two further zooms near the boundary of Fig~\ref{fig:1.4}}
\end{figure}

In contrast to these examples, Figures~\ref{fig:1.4} and~\ref{fig:1.5} display a rather different phenomenon, namely
more and more bands emanating from the red (scattering) region as one zooms into small and smaller scales.
To be more specific, the image on the right-hand side
 of Figure~\ref{fig:1.4} as well as those  in Figure~\ref{fig:1.5} show three zooms with successively higher resolutions
  near the round tip 
 of the left-hand image of Figure~\ref{fig:1.4}.  
 
 With some patience the reader will not only be able to find the filaments that are visible in one scale
 again at the next finer scale, but also between those and the solid red region a very thin new filament emerging.
 Note that the green color in the three zooms signifies the ``indecisive" points where the algorithm was not able
 to decide  between blowup and global existence.  It is this indecisive region that gets resolved at the next finer
 scale, with a new filament emerging from it. 

 The zooms were computed  on the Kraken
 supercomputer at the NICS, which is part of the TeraGrid.  
 The coarse appearance of the third zoom, shown on the right-hand side of Figure~\ref{fig:1.4} is due to 
the limitations on the computing time available to the authors. 

At this point, it is completely unclear whether or not infinitely many of these ``bands'' emerge out of
the scattering region. Of course, this would imply that $\p\S_+$ is {\em not} a smooth manifold --- perhaps a
somewhat counter-intuitive assertion. 

A very visible common feature of all figures of this section (excluding of course the zooms)
 are the thin ``spikes'' that emanate from the bulk of the red region (the forward scattering set $\p\S_{+}$). 
Even though current theory cannot touch such phenomena as these ``spikes'', we now offer some --- admittedly
daring --- explanations as to why they might appear. First, let us assume the {\em scattering hypothesis} which
states that solutions  starting from data that fall {\em exactly} onto the boundary exist globally in forward time
and scatter to either $Q$ or $-Q$.  In addition to the proof of this hypothesis in~\cite{NakS1} for energies
which are at most slightly larger than that of the ground state, circumstantial evidence 
gathered by means of the bisection method described in the previous section
lends some credence to this conjecture. However, if there are indeed infinitely many filaments in Figures~\ref{fig:1.3}, \ref{fig:1.4} then
it is very difficult to guess  what the dynamical implications might be; in particular, it is hard to conjecture (or simulate for
that matter) what kind of  asymptotic states might be associated with such a complicated boundary. 

Be it as it may, let us assume the ``scattering hypothesis''. Then the boundary splits into two open sets which scatter to
either $Q$ or $-Q$ under the nonlinear flow (for convenience, let us call these {\em positive} and {\em negative} boundaries, respectively). 
The openness of these sets is another leap of faith, but a natural one: it simply expresses that the scattering properties should 
be stable {\em relative to small perturbations inside the boundary}. 
If    this is so then these sets will not meet up smoothly but rather avoid each other. 
If a two-dimensional section does indeed display smooth boundary curves as in the numerical pictures, then 
the positive and negative  curves cannot approach each other within a compact set. Therefore, they have to move off
to infinity. Note that this does not exclude that two boundary curves form a ``tube'' of some finite thickness, rather than a spike.
Such a tube, however, is extremely unlikely since it would insinuate something like a uniform stability property at larger and larger energies. 
 
In this context, it would be very instructive to consider the Schr\"odinger equation since the soliton manifold in that case is
connected (whereas for the NLKG equation it is discrete: $\{\pm Q\}$). 
If the previous explanation has any merit, then one might expect that spikes do not form in the Schr\"odinger case. 
Some preliminary numerical experiments by the authors appear to confirm that. More precisely, two-dimensional sections of
the  scattering region for the cubic radial focusing NLS equation in $\R^{3}$ 
appear as round and smooth shapes --- at least in the very limited number of examples
the authors computed for NLS. 

We remark that the ``spikes'' should not be confused with the one-dimensional stable and unstable manifolds emanating from $\pm Q$ (or rather,
neighborhoods thereof). Indeed, those manifolds   would
be only visible in very carefully chosen curved two-dimensional sections  passing through $\pm Q$, and finding such a section would be a zero probability
event in any reasonable sense.  But this contradicts the fact that the spikes appear in each of  the numerically computed images (but not in the zooms, of course)  that the authors produced.  

Needless to say, it is essentially impossible  to hit the boundary exactly 
just by sampling a grid of data as we do to produce the figures of this section. 
Nevertheless, one can   
test the aforementioned scattering  hypothesis by means of the ``bisection method'' described in Section~\ref{sec:code}, see also Bizo\'n et al.~\cite{Biz}, \cite{BCS}. 
This refers to taking two data pairs, one from the red region, and another from the white region which are in close
proximity of each other. Then one numerically computes the solution starting from the midpoint between these two points, tests for blowup/dispersion,
and continues bisecting the interval whose endpoints exhibit opposite behaviors. 
In this way one obtains data that have a clear tendency of approaching either $Q$ or $-Q$ for large times, albeit with
possible oscillations about these two states, followed by eventual blowup or dispersion. The aforementioned oscillations are referred to
as ``quasinormal ringing'' in the physics literature, see for example Bizo\'n et al. 

To the authors, this bisection process also provides a strong check of the accuracy and validity of the numerically computed region $\S_{+}$. 
Indeed, if the initial data pair for the bisection process does not lie close to, and on opposite sides of, the {\em actual} $\p\S_{+}$ --- which is 
what the figures claim --- then the bisection process would not converge in the described fashion.

\begin{figure}[ht]
\begin{center}
 \includegraphics[width=0.495\textwidth]{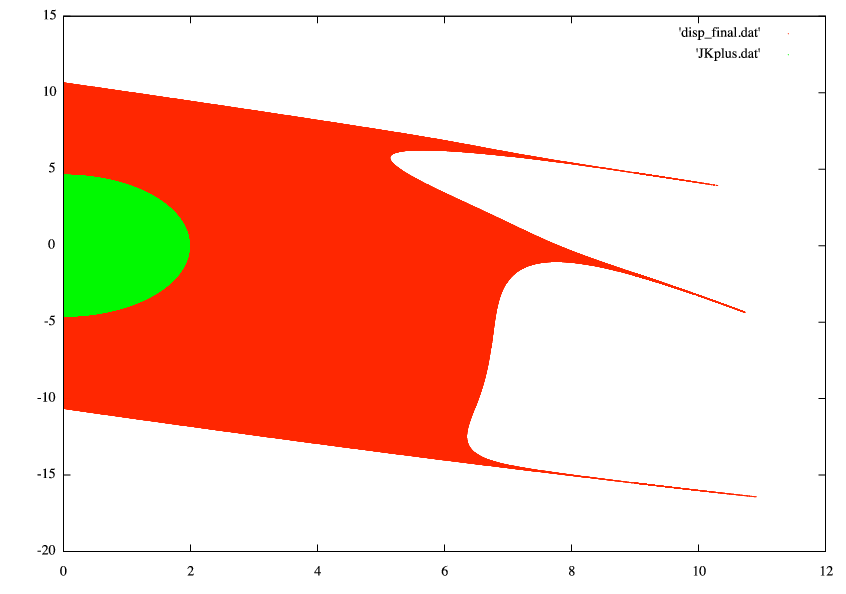}
 \includegraphics[width=0.495\textwidth]{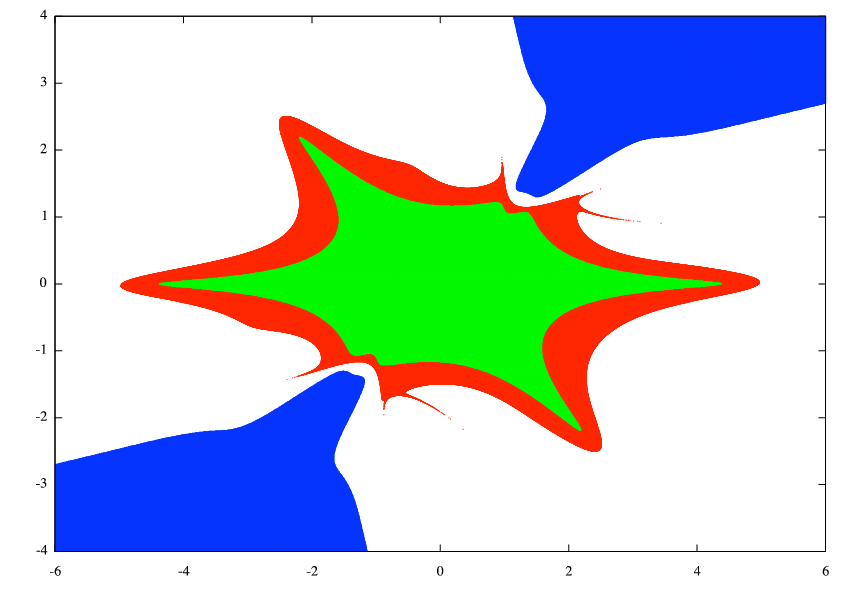}
\end{center}
 \caption{ \label{fig:1.6} $\D_{2} ( r\sin(6r) e^{-\f32\lan r\ran}, \cos(6r) e^{-(r^{2}-\f14)^{2}}  )$ and the curved section from~\eqref{curve1}}
\end{figure}

Figure~\ref{fig:1.6} shows two more sections with analytic, decaying data. The left-hand image is a planar
section as defined above, and it provides (after the reflection about the origin) an example with {\em six spikes}
rather than two as in the previous examples. A priori, there seems to be no reason to assume that the number of spikes
could not be arbitrarily large.   

The image on the right-hand side of Figure~\ref{fig:1.5} is different from all those previously considered in so far
as {\em it is not planar}. More specifically, it shows a section produced by 
\EQ{\label{curve1}
(f,g) =\big( (A+B)B\cos(2(A-B)r) \exp(-\f{r^{2}}{\lan r \rangle}),     (A-B)\exp(-r^{2})   \big), 
}
with $|A|\le 6, \; |B|\le 4$. 
%

\begin{figure}[ht]
\begin{center}
 \includegraphics[width=0.495\textwidth]{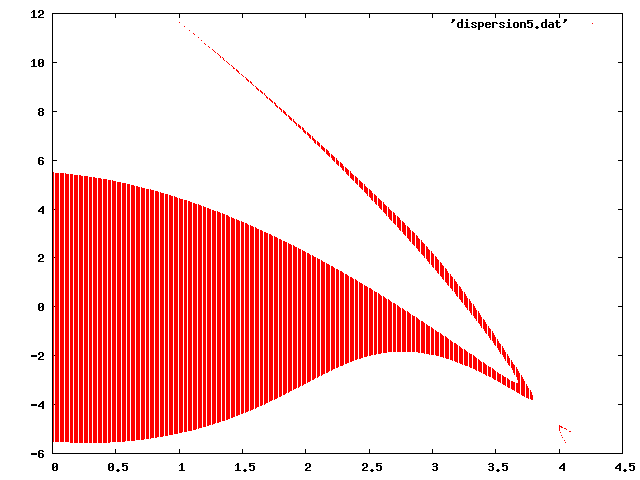}
 \includegraphics[width=0.495\textwidth]{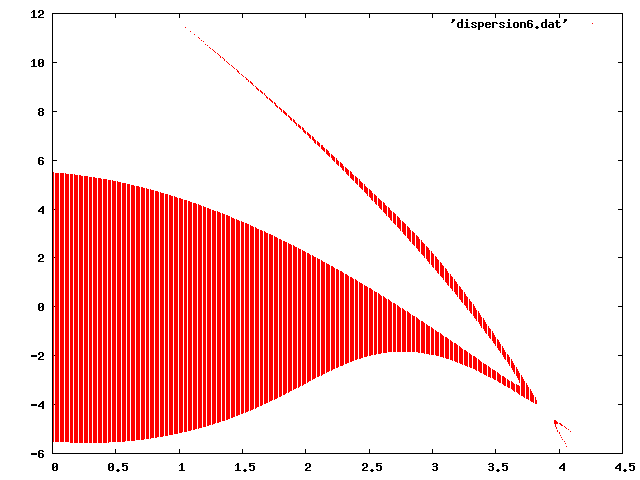}
  \includegraphics[width=0.495\textwidth]{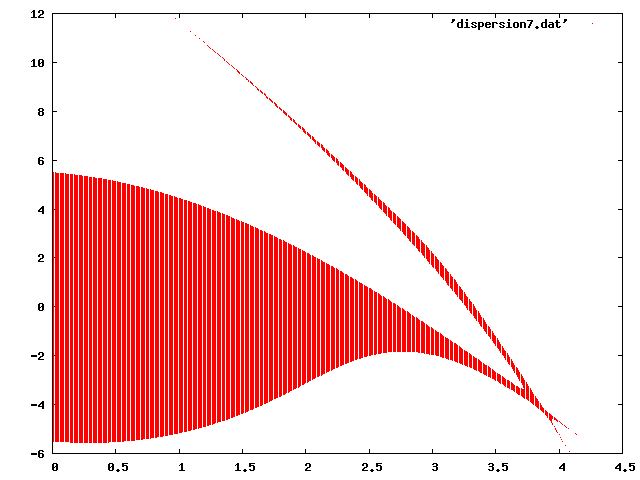}
   \includegraphics[width=0.495\textwidth]{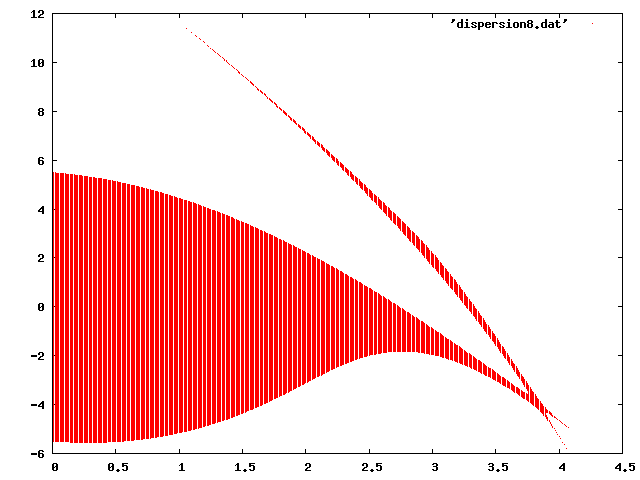}
 \includegraphics[width=0.495\textwidth]{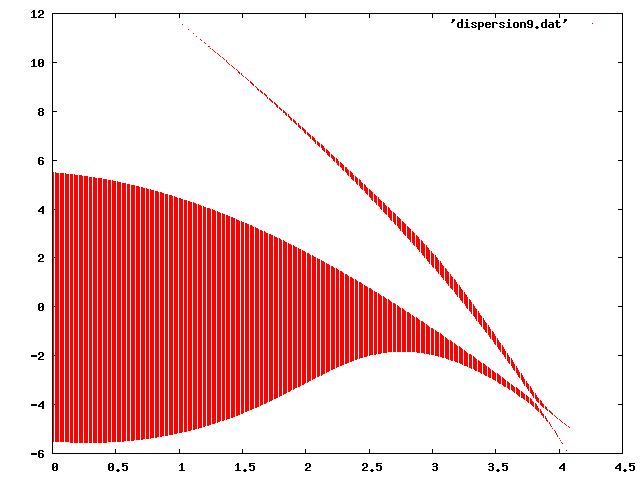}
  \includegraphics[width=0.495\textwidth]{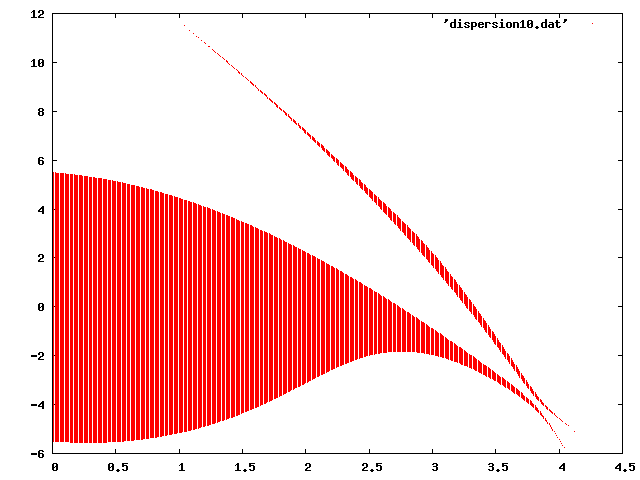}
\end{center}
 \caption{ \label{fig:1.7}  A series of sections as in \eqref{3param} with increasing $C$}
\end{figure}

\begin{figure}[ht]
\begin{center}
 \includegraphics[width=0.32\textwidth]{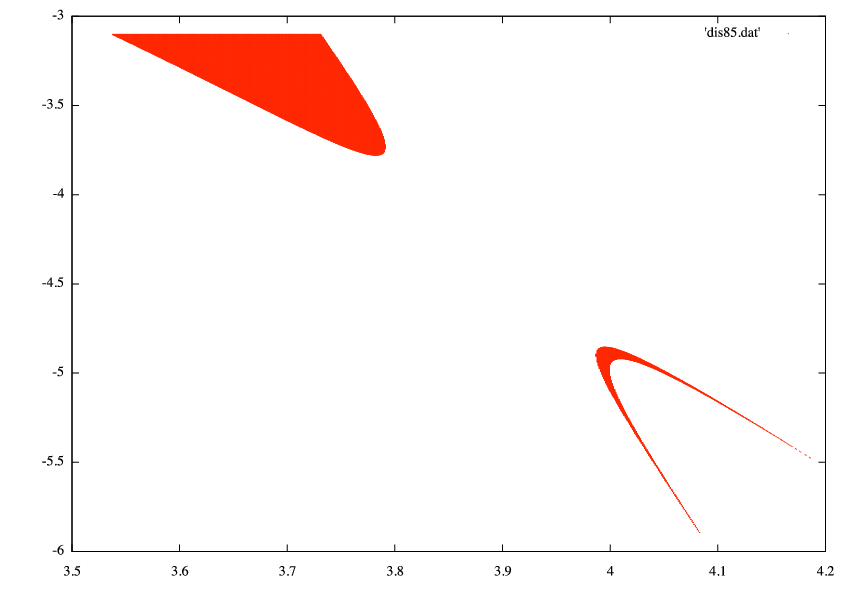}
 \includegraphics[width=0.32\textwidth]{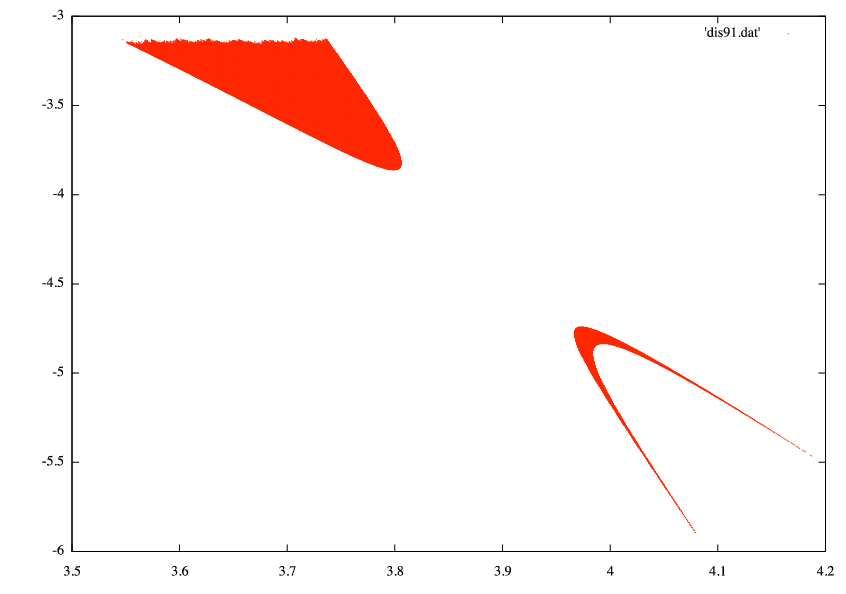}
  \includegraphics[width=0.32\textwidth]{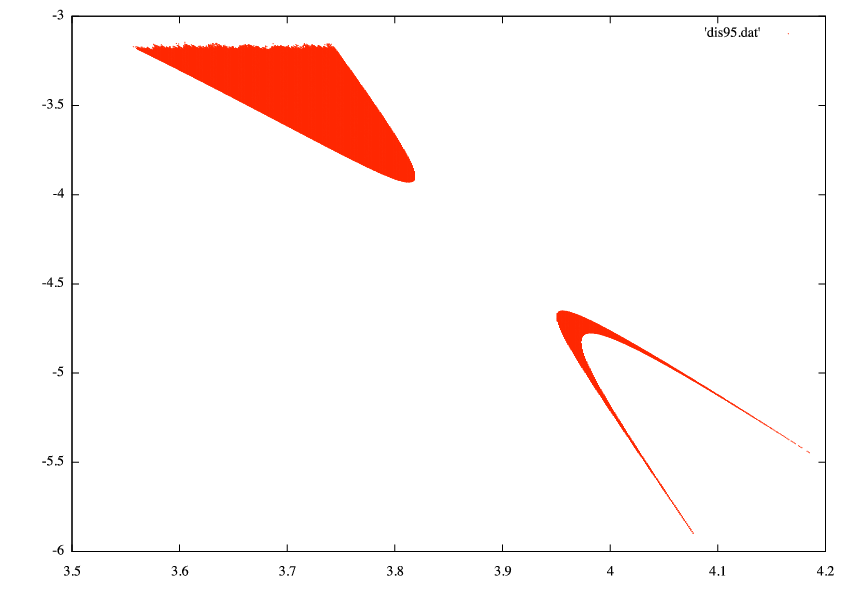}
   \includegraphics[width=0.32\textwidth]{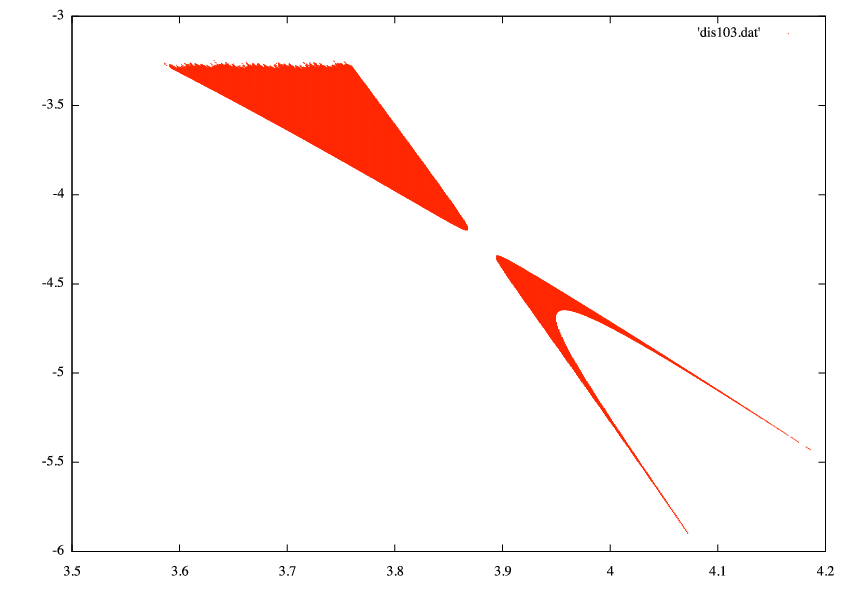}
 \includegraphics[width=0.32\textwidth]{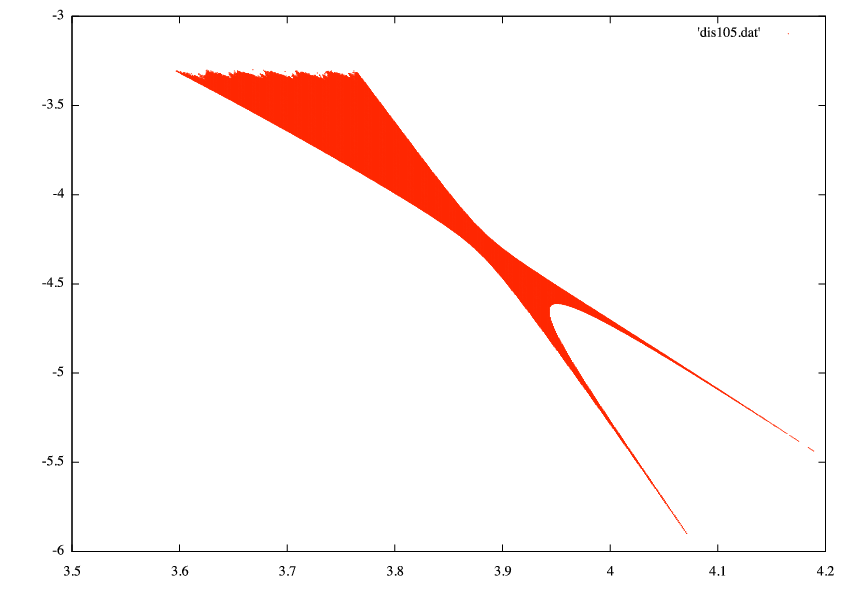}
  \includegraphics[width=0.32\textwidth]{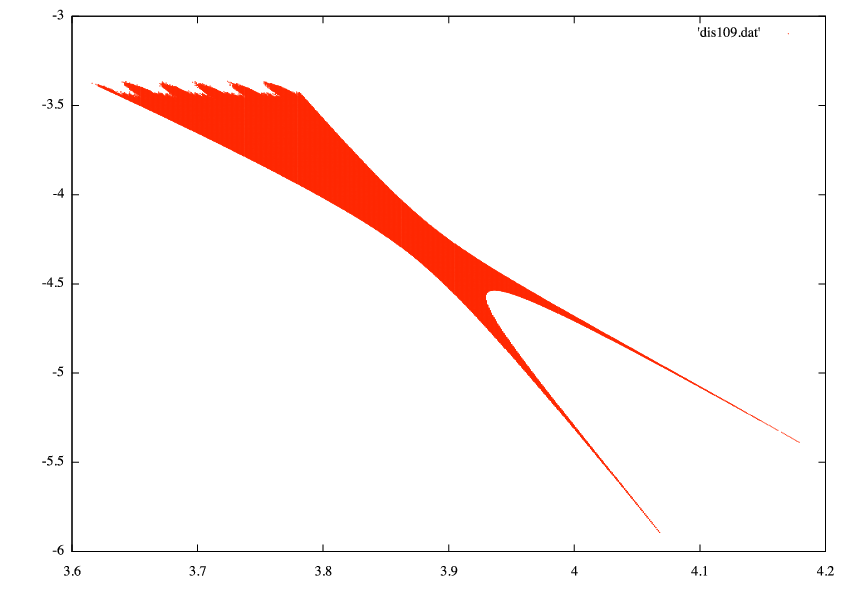}
  \includegraphics[width=0.32\textwidth]{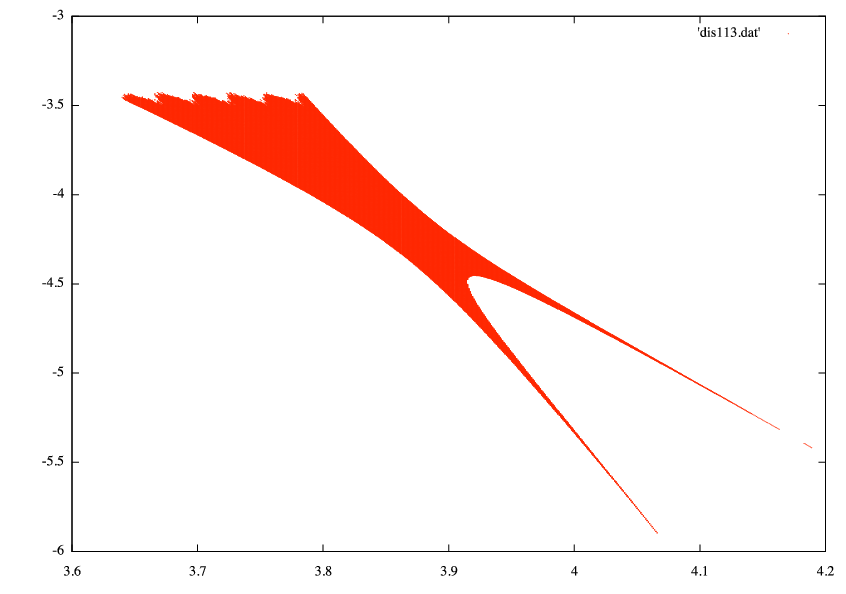}
   \includegraphics[width=0.32\textwidth]{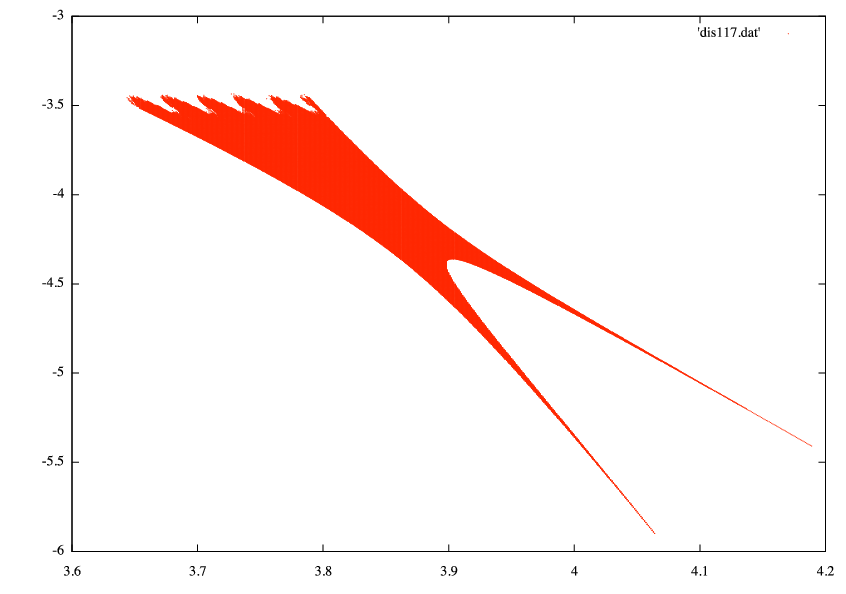}
 \includegraphics[width=0.32\textwidth]{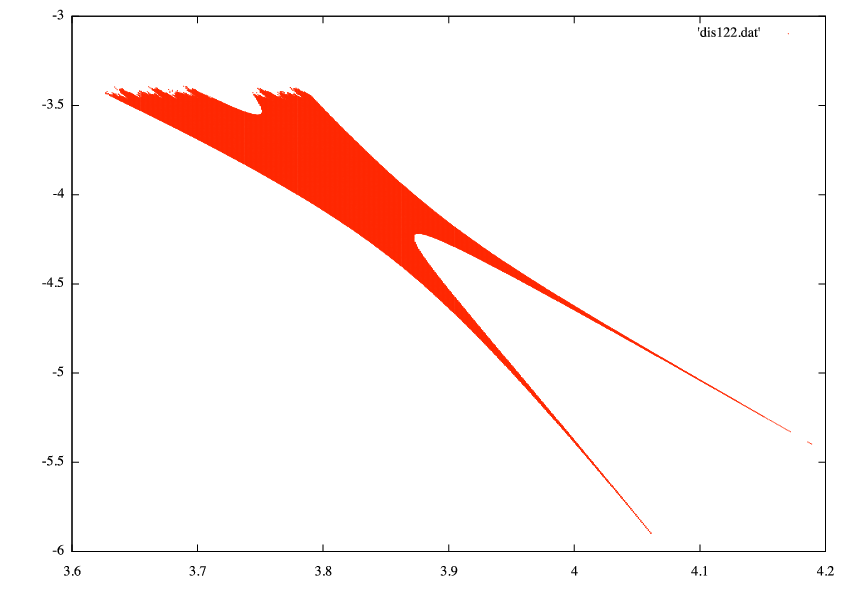}
   \includegraphics[width=0.32\textwidth]{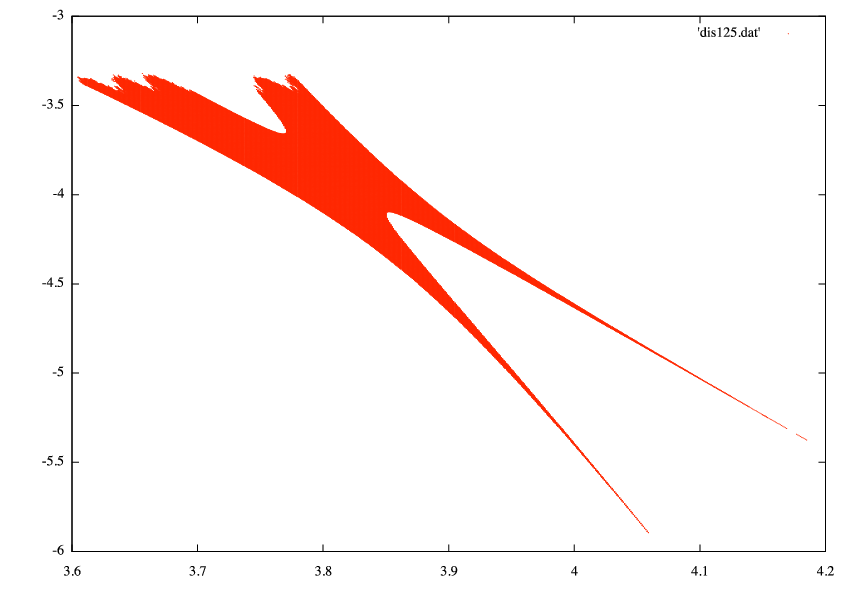}
  \includegraphics[width=0.32\textwidth]{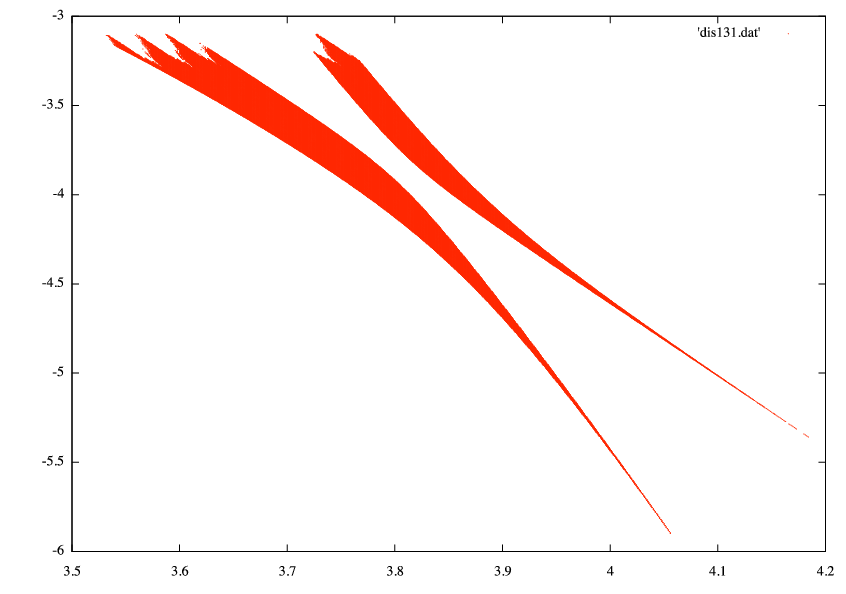}
  \includegraphics[width=0.32\textwidth]{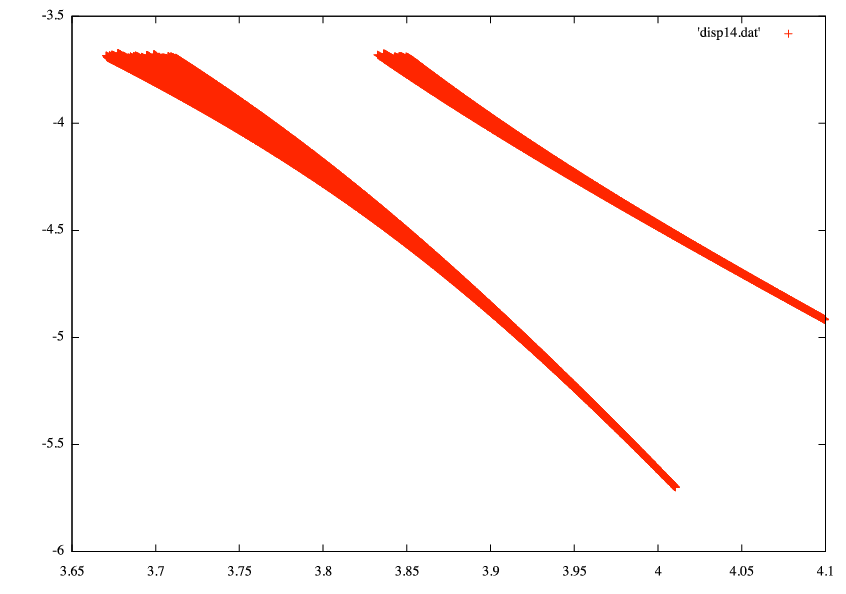}
\end{center}
 \caption{ \label{fig:1.8}  A series of higher-resolution images for \eqref{3param} with increasing $C$}
\end{figure}

The final series of numerical experiments which we present here concerns {\em three-parameter families}.
This refers to data of the form $(Af(C,r), B g(C,r))$, say,  where $A,B,C$ are parameters which are chosen in some finite range.
For fixed $C$ these data generate a two-dimensional section as shown in the previous images, but as $C$ varies
one obtains a whole family of sections which move ``continuously'' (at least, if $C$ changes only by very little). 
The specific data leading to  Figures~\ref{fig:1.7}, \ref{fig:1.8}  are
\EQ{\label{3param}
\D_{2} \big( r^{2}\cos(Cr) e^{-r^{2}}, r^{3}\sin(6r) \exp(-2(r^{2}-1)^{2} \big)
}
where $C$ is in the vicinity of $9.72$, $9.73$ and up. Figure~\ref{fig:1.7} shows a series of such two-dimensional slices
obtained by letting $C$ increase through six values (the white vertical lines visible in the red region are an artifact caused by  a slightly
coarser resolution in the horizontal direction). What can be seen in these six images, starting with the upper
left, is first the merging of two disjoint sets, followed by the pulling apart of two sets which then become disjoint.
Moreover, the direction in which the sets pull apart is roughly orthogonal to that in which the merging takes place. 

The images in Figure~\ref{fig:1.7} were computed on the MacPro and use a fairly low resolution
in order to keep the computation times reasonable. A much finer resolution was obtained on the supercomputer
Kraken at the NICS, leading to the sequence of images shown in Figure~\ref{fig:1.8}. They are zoom-ins into
the area where the pinching takes place as shown in Figure~\ref{fig:1.7}, and depict the process of merging
and pulling apart more clearly. The uneven appearance of some of the upper boundaries in the images of  Figures~\ref{fig:1.8} is
due to the fact that in most computations the so-called wall-time was reached which then leads to termination. 

Interpreting these images is a delicate issue, as they are two-dimensional sections of an infinite-dimensional object, namely $\S_{+}$.
Note carefully that they certainly {\em do not} imply singularity formation of the actual boundary $\p\S_{+}$ as can be seen from
the example of the one-sheeted hyperboloid in $\R^{3}$, say with rotation axis given by the $x_{3}$-coordinate direction.
Then as the planes $x_{1}=\const$ slice through this hyperboloid they will produce an image given by $x_{2}^{2}-x_{3}^{2}\le0$ at
the moment of tangency to the waist line.  This is what seems to happen between the fourth and fifth image of Figure~\ref{fig:1.8}.

\section{Concluding remarks}
\label{sec:coda}

The authors hope that this note will encourage further numerical, or analytical,  investigations of the set $\S_{+}$, both for the NLKG equation in various dimensions and with various powers in the nonlinearity, 
as well as other nonlinear dispersive equations such as the nonlinear Schr\"odinger equation.  It would be highly desirable to employ different, and possibly more sophisticated, numerical algorithms. 
For example, one could try  the symplectic difference schemes developed by Duncan~\cite{Duncan} and Reich~\cite{Reich}, or a different approach altogether. 

The most important problem arising in this context is to decide whether or not the boundary of the forward scattering region can indeed become singular at larger
energies. And if so, in what ways the boundary may degenerate. As discussed in the previous section, the only evidence we have found that
points towards the possibility of some sort of singular behavior is the emergence of more and more filaments from the scattering region. 
The next step would then be to see if any possible degeneration of the boundary is reflected by anomalies in the {\em actual wave dynamics}. 

Upon request the authors will email all numerical data which they have produced in the course of this work (which is considerably larger than that presented here). 
However, hopefully  other researchers will carry out their experiments  with data and numerical solvers of their own choosing. 
It seems most important to conduct a systematic and exhaustive study of the boundary of $\S_{+}$ by means of the bisection method
which we discussed above in order to test the ``scattering hypothesis''. In the process one can also 
test the accuracy of the numerically computed
 two-dimensional sections, particularly if a different solver is employed for the bisection than the one used for the computation of 
 the two-dimensional slices of~$\S_{+}$.  
 
 More work needs to be done on  three-parameter families.
 In particular, it would be of interest to study how the ``bubbles'' in Figure~\ref{fig:1.3} evolve as a third parameter is being changed.  
 Furthermore, one could use such a three parameter family to obtain a two-dimensional section of\footnote{We thank Kenji Nakanishi for this suggestion.} 
 {\em constant energy}. And thirdly, it would be very appealing to produce three-dimensional images
 of the intersections of $\S_{+}$ with a cube in $\R^{3}$ (this requires substantial computational resources).

The issue here is not so much to produce more pictures, but to pinpoint various phenomena relating to the topology
of $\p\S_{+}$ and {\em then to connect them to the actual wave dynamics}. An example of this is the aforementioned {\em scattering
hypothesis} but one can also ask about the dynamical implications (if any) of the filaments as they appear in 
Figures~\ref{fig:1.4} and~\ref{fig:1.5}. 

Finally, it would be interesting --- as well as challenging --- to consider the same questions that we pursue here for other
equations, such as   the one-dimensional NLKG, the Schr\"odinger, as well as the energy critical wave equations (i.e., with $u^{5}$ nonlinearity)  in $\R^{3}$. The one-dimensional equation
 has considerably less dispersion than the three-dimensional one, so it is perhaps somewhat more delicate to rely on the numerics in that case.
 In this context, one can also  investigate the dependence on  the power of the nonlinearity, see Bizon~\cite{BCS} for
 such considerations in connection with the aforementioned ``quasi-periodic ringing'' and the emergence of a threshold resonance at power $p=3$. 
 It is natural to consider even data as in~\cite{KNS2}, but general data are also important, of course, and considerably more challenging numerically.
 As for the Schr\"odinger and critical wave equations, these are numerically (as well as analytically)  harder and the authors have done only
 preliminary exploratory work on them as far as the numerics is concerned.

\end{document}